\documentclass[aps,reprint,letter,nofootinbib]{revtex4-1}
\usepackage{amssymb,amsmath}
\usepackage{graphicx}
\usepackage{bm}
\usepackage{hyperref}

\def\noterule{\bigskip\hrule height 1.2pt\medskip}

\newcommand{\En}{{\cal{E}}}

\newcommand\Tr{T}
\newcommand{\GR}{G}

\newcommand{\Trf}{\Tr_f}

\newcommand{\levy}{L{\'e}vy}

\newcommand{\Ga}{\alpha}
\newcommand{\Gb}{\beta}
\newcommand{\Gd}{\delta}

\newcommand{\Gf}{\phi}
\newcommand{\Gg}{\gamma}
\newcommand{\Gk}{\kappa}
\newcommand{\Gl}{\lambda}
\newcommand{\GL}{\Lambda}
\newcommand{\Gm}{\mu}

\newcommand{\GP}{\Phi}

\newcommand{\dd}{{\hbox{\rm d}}}

\newcommand{\Eq}[1]{Eq.~(\ref{eq:#1})}

\newcommand{\ovr}[2]{{{#1}\over{#2}}}

\newcommand{\prt}{\partial}
\newcommand{\povr}[2]{\ovr{\prt #1}{\prt #2}}

\begin{document} 

\title{A simple derivation and classification of common probability distributions based on information symmetry and measurement scale}

\author{Steven A.~Frank}
\affiliation{Department of Ecology and Evolutionary Biology,
  University of California, Irvine, CA 92697-2525, USA}
\affiliation{Santa Fe Institute, 1399 Hyde Park Road, Santa Fe, NM
87501, USA}
\author{Eric Smith}
\affiliation{Santa Fe Institute, 1399 Hyde Park Road, Santa Fe, NM
87501, USA}

\date{\today}

\begin{abstract}
Commonly observed patterns typically follow a few distinct families of probability distributions. Over one hundred years ago, Karl Pearson provided a systematic derivation and classification of the common continuous distributions.  His approach was phenomenological: a differential equation that generated common distributions without any underlying conceptual basis for why common distributions have particular forms and what explains the familial relations.  Pearson's system and its descendants remain the most popular systematic classification of probability distributions.  Here, we unify the disparate forms of common distributions into a single system based on two meaningful and justifiable propositions.  First, distributions follow maximum entropy subject to constraints, where maximum entropy is equivalent to minimum information.  Second, different problems associate magnitude to information in different ways, an association we describe in terms of the relation between information invariance and measurement scale.  Our framework relates the different continuous probability distributions through the variations in measurement scale that change each family of maximum entropy distributions into a distinct family.
\end{abstract}

\maketitle

\section{Introduction}

Commonly observed patterns follow a few families of probability distributions.  For example, Gaussian patterns often arise from measures of height or weight, and gamma patterns often arise from measures of waiting times.  These common patterns lead to two questions.  How are the different families of distributions related?  Why are there so few families, when the possible patterns are essentially infinite?  

These questions are important, because one can hardly begin to study nature without some sense of the fundamental contours of pattern and why those contours arise.  For example, no one observing a Gaussian distribution of weights in a population would feel a need to give a special explanation for that pattern.  The central limit theorem tells us that a Gaussian distribution is a natural and widely expected pattern that arises from measuring aggregates in a certain way.  

With other common patterns, such as power laws, the current standard of interpretation is much more variable. That variability arises because we do not have a comprehensive theory of how measurement and information shape the commonly observed patterns.  Without a clear notion of what is expected in different situations, common and relatively uninformative patterns frequently motivate unnecessarily complex explanations, and surprising and informative patterns may be overlooked \cite{frank09the-common}. 

Currently, the differences between families of common probability distributions often seem arbitrary.  Thus, little understanding exists with regard to how changes in process or in methods of observation may cause observed pattern to change from one common form into another.  

We argue that measurement, described by the relation between magnitude and information, unifies the distinct families of common probability distributions.  Variations in measurement scale may, for example, arise from varying precision in observations at different magnitudes or from the way that information is lost when measurements are made on aggregates.  Our unified explanation of the different commonly observed distributions in terms of measurement points the way to a deeper understanding of the relations between pattern and process. 

We develop the role of measurement through maximum entropy expressions for probability distributions.  We first note that all probability distributions can be expressed by maximization of entropy subject to constraint. Maximization of entropy is equivalent to minimizing total information while retaining all the particular information known to constrain underlying pattern \cite{jaynes57information,jaynes57informationII,jaynes03probability}. To obtain a probability distribution of a given form, one simply chooses the informational constraints such that maximization of entropy yields the desired distribution.  However, constraints chosen to match a particular distribution only describe the sufficient information for that distribution.  To obtain deeper insight into the causes of particular distributions and each distribution's position among related families of distributions, we derive the related forms of constraints through variations in measurement scale.  

Measurement scale expresses information through the invariant transformations of measurements that leave the form of the associated probability distribution unchanged \cite{frank10measurement}.  Each problem has a characteristic form of information invariance and symmetry that sets the measurement scale \cite{hand04measurement,luce08measurement,narens08meaningfulness} and the most likely probability distribution associated with that particular scale \cite{frank10measurement}.  We show that measurement scales  and the symmetries of information invariances form a natural hierarchy that generates the common families of probability distributions.  We use \textit{invariance\/} and \textit{symmetry\/} interchangeably, in the sense that symmetry arises when an invariant transformation leaves an object unchanged \cite{weyl52symmetry}.

The measurement hierarchy arises from two processes.  First, we express the forms of information invariance and measurement scale through a continuous group of transformations, showing the relations between different types of information invariance.  Second, the types of aggregation and measurement that minimize information and maximize entropy fall into two classes, each class setting a different basis for information invariance and measurement scale. 

The two types of aggregation correspond to the two major families of stable distributions that generalize the process leading to the central limit theorem: the \levy\ family that includes the Gaussian distribution as a special case, and the Fisher-Tippett family of extreme value distributions.  By expressing measurement scale in a general way, we obtain a wider interpretation of the families of stable distributions and a broader classification of the common distributions.

Our derivation of probability distributions and their familial relations supersedes the Pearson and similar classifications of continuous distributions \cite{johnson94continuous}.  Our system derives from a natural description of varying information in measurements under different conditions \cite{frank10measurement}, whereas the Pearson and related systems derive from phenomenological descriptions that generate distributions without clear grounding in fundamental principles such as measurement and information. 

Some recent systems of probability distributions, such as the unification by Morris \cite{morris09unifying,morris82natural}, provide great insight into the relations between families of distributions.  However, Morris's system and other common classifications do not derive from what we regard as fundamental principles, instead arising from descriptions of structural similarities among distributions.  We provide a detailed analysis of Morris's system in relation to ours in Appendix C.

We favor our system because it derives the relations between distributions from fundamental principles, such as maximum entropy and the invariances that define measurement scale.  Although the notion of what is fundamental will certainly attract controversy, our favored principles of entropy, symmetries defined by invariances, and measurement scale certainly deserve consideration.  Our purpose is to show what one can accomplish by starting solely with these principles.

\section{Maximum entropy and measurement scale}

This section reviews our prior work on the roles of information invariance and measurement scale in setting observed pattern \cite{frank10measurement}.  The following sections extend this prior work by expressing measurement in terms of the scale of aggregation and the continuous group transformations of information invariance.

\subsection{Maximum entropy}

The method of maximum entropy defines the most likely probability distribution as the distribution that maximizes a measure of entropy (randomness) subject to various information constraints \cite{jaynes03probability}.  We write the quantity to be maximized as
\begin{equation}\label{eq:maxEnt}
	\GP = \En - \Gk C_0 - \sum_{i=1}^n\Gl_iC_i,
\end{equation}
where $\En$ measures entropy, the $C_i$ are the constraints to be satisfied, and $\Gk$ and the $\Gl_i$ are the Lagrange multipliers to be found by satisfying the constraints.  Let $C_0=\int p_y\dd y -1$ be the constraint that the probabilities must total one, where $p_y$ is the probability distribution function of $y$.  The other constraints are usually written as $C_i= \int p_yf_i(y)\dd y -\bar{f}_i$, where the $f_i(y)$ are various transformed measurements of $y$, and the overbar denotes mean value. A mean value is either the average of some function applied to each of a sample of observed values, or an a priori assumption about the average value of some function with respect to a candidate set of probability laws. If $f_i(y)=y^i$, then $\bar{f}_i$ are the moments of the distribution---either the moments estimated from observations or a priori values of the moments set by assumption.  The moments are often regarded as ``standard'' constraints, although from a mathematical point of view, any properly formed constraint can be used.  

Here, we confine ourselves to a single constraint of measurement. We express that constraint with a more general notation, $C_1= \int p_y\Tr(f_y)\dd y -\bar{\Tr}_f$, where $f_y\equiv f(y)$, and $\Tr(f_y)\equiv\Trf$ is a transformation of $f_y$.  We could, of course, express the constraining function for $y$ directly through $f_y$.  However, we wish to distinguish between an initial function $f_y$ that can be regarded as a standard measurement, in any sense in which one chooses to interpret the meaning of standard, and a transformation of standard measurements denoted by $\Trf$ that arises from information about the measurement scale.  

The maximum entropy distribution is obtained by solving the set of equations
\begin{equation}\label{eq:maxEntSoln}
	\povr{\GP}{p_y} = \povr{\En}{p_y} - \Gk - 
	\Gl\Trf=0,
\end{equation}
where one checks the candidate solution for a maximum and obtains $\Gk$ and $\Gl$ by satisfying the constraint on total probability and the constraint on $\bar{\Tr}_f$. We assume that we can treat the entropy measures and the maximization procedure by the continuous limit of the discrete case.

In the standard approach, we define entropy by extension of Shannon information 
\begin{equation}\label{eq:shannonDef}
	\En=-\int p_y\log\left(\ovr{p_y}{m_y}\right)\dd y,
\end{equation}
in which this expression may be called Jaynes's differential entropy \cite{jaynes03probability}, which is equivalent in form to the continuous expression of relative entropy or the Kullback-Leibler divergence \cite{cover06elements}. Here, we will interpret $m_y$ by information invariance and measurement scale as discussed below.   
With these definitions, the solution of \Eq{maxEntSoln} is
\begin{equation}\label{eq:shannonSoln}
	p_y \propto m_ye^{- \Gl\Trf},
\end{equation}
where $\Gl$ satisfies the constraint $C_1$, and the proportionality is adjusted so that the total probability is one by choosing the parameter $\Gk$ to satisfy the constraint $C_0$.

\subsection{Information invariance and measurement scale}\label{sectionInvar}

Maximum entropy must capture all of the available information about a particular problem.  One form of information concerns transformations to the measurement scale that leave the most likely probability distribution unchanged \cite{jaynes03probability,frank09the-common,frank10measurement}.  Here, it is important to distinguish between measurements and measurement scale.  In our notation, we start with measurements, $f_y$, made on the measurement scale $y$.  For example, we may have measures of squared deviations about zero, $f_y=y^2$, with respect to the measurement scale $y$, such that $\bar{f}_y$ is the second moment of the measurements with respect to the underlying measurement scale.  

Suppose that we obtain the same information about the underlying probability distribution from measurements of $f_y$ or transformed measurements, $\GR(f_y)$.  Put another way, if one has access only to measurements $\GR(f_y)$, one has the same information that would be obtained if the measurements were reported as $f_y$.  We say that the measurements $f_y$ and $\GR(f_y)$ are equivalent with respect to information, or that the transformation $f_y \rightarrow \GR(f_y)$ is an information invariance that describes a symmetry of the measurement scale.

To capture this information invariance in maximum entropy, we must express our measurements so that 
\begin{equation}\label{eq:transDef}
  \Tr(f_y) = \Gd + \Gf\Tr[\GR(f_y)]
\end{equation}
for some arbitrary constants $\Gd$ and $\Gf$ \cite{frank10measurement}.  Putting this definition of $\Tr(f_y)\equiv\Trf$ into \Eq{shannonSoln} shows that we get the same maximum entropy solution whether we use the observations $f_y$ or the transformed observations, $\GR(f_y)$, because the $\Gk$ and $\Gl$ constants will adjust to the constants $\Gd$ and $\Gf$ so that the distribution remains unchanged.

\section{Deriving probability distributions}

The prior section established two key steps.  First, maximum entropy probability distributions have the form given in \Eq{shannonSoln} as $p_y \propto m_ye^{- \Gl\Trf}$.  Second, the expression of $\Tr(f_y)$ for each problem comes from the particular information invariance $\GR(f_y)$ associated with that particular problem.  To derive specific probability distributions, we must pass three further steps, which we develop in the following sections.

First, we turn the abstract notions of information invariance and measurement scale into specific expressions for the measurement scale function, $\Tr(f_y)$.  We accomplish this by developing the continuous group transformations for information invariance.  Those continuous transformations provide an abstract hierarchy of forms for probability distributions based on the scale factor, $m_y$, the specific measured attribute, $f_y$, and how the information and precision of measurements change with magnitude expressed by the measurement scale $\Tr(f_y)$.

Second, we define $m_y$ as the relation between the scale of information invariance and the scale on which we express probability.  To use the maximization of entropy and the associated minimization of information, we must relate the information invariance of measurement to the scale on which underlying processes dissipate information.  We consider alternative interpretations of scale that may be associated with the dissipation of information by aggregation of random perturbations and by measurements of extreme values. We also consider measurements on a scale that differs from the basis for dissipation of information.

Third, we consider how to interpret $f_y$, which is the value used to describe the informational constraint in relation to the measurement scale $\Tr(f_y)$, leading to the constraint $\bar{\Tr}_f$.  We discuss $f_y$ as a reduction in the dimensionality of information to a single sufficient dimension. That sufficient dimension sets the form of probability under the various processes of information dissipation that lead to the common probability distributions.

\section{Continuous group transformations of measurement}

The transformation in \Eq{transDef} sets the relation between information invariance and measurement scale.  However, that expression does not show in a simple way the relations between information and measurement.

To understand commonly observed patterns in relation to the families of probability distributions, it is helpful to express in a general way the underlying symmetry that determines information invariance and measurement scale.  From that underlying symmetry, we may see more clearly the associated relations between the forms of probability distributions.  

\subsection{The affine structure}

The relation between information invariance and measurement scale in \Eq{transDef} arises directly from the form of maximum entropy solutions in \Eq{shannonSoln}, in which probability distributions are exponentials of the transformed constraint measures, $\Tr_f$.  In particular, the probability distribution associated with a constraint is invariant to an additive shift of the constraint and a multiplicative change in the scale of the constraint, given by the parameters $\Gd$ and $\Gf$ in \Eq{transDef}.  It is that symmetry in the affine structure of invariant transformation that ultimately sets the underlying relations between information, measurement, and the familial forms of the common probability distributions.

To understand the affine structure of the invariant transformation in \Eq{transDef} more clearly, we can express that invariant transformation as a continuous operator.  First, rearrange \Eq{transDef} as an equivalent expression
\begin{equation}\label{eq:transDefInv}
  \Tr[\GR(f_y)] = a + b\Tr(f_y)
\end{equation}
with new parameters $a$ and $b$ that are easily calculated from \Eq{transDef}.  We show in Appendix A that we can express the same information invariance of $\GR(f_y)$ by the differential operator defined as
\begin{equation}
  v_w =  
  \left( \alpha + \beta \Tr \right)
  \frac{\dd}{\dd\Tr}  
\label{eq:v_alpha_beta_form}
\end{equation}
that can be applied to $\Tr$ as
\begin{equation}
  v_w \! \left( \Tr \right) = 
  \alpha + \beta\Tr . 
\label{eq:v_alpha_beta_onT}
\end{equation}
Recursive application of $v_w$ preserves the affine structure and so keeps the successive transformations within the group of admissible invariance relations.

We can express $v_w$ as
\begin{equation}
  v_w = 
  \frac{\dd}{\dd w} , 
\label{eq:v_w_form}
\end{equation}
where $w\equiv w (f_y)$ is some function of $f_y$.  We then have a differential equation for $\Tr$ as
\begin{equation}
  \frac{\dd\Tr}{\dd w} - 
  \beta \Tr = 
  \alpha , 
\label{eq:T_resp_genform}
\end{equation}
which has solutions of the general form 
\begin{equation}
  \Tr \! \left( f_y\right) = 
  \Tr_0
  e^{\beta w} + 
  \frac{\alpha}{\beta}
  \left( 
    e^{\beta w} - 1 
  \right) , 
\label{eq:T_gen_soln}
\end{equation}
which as $\beta \rightarrow 0$ goes to $\Tr(f_y) \rightarrow \Tr_0 + \alpha w$.  \Eq{T_gen_soln} gives the most general class of measurement functions, $\Tr(f_y)$, for which the associated transformations generated by $v_w$ preserve information invariance.

The operator $v_w$ can be applied repeatedly, creating a recursively generated sequence of deformations that all satisfy the fundamental relation between deformations of measurement and information invariance.  By thinking of $w(f_y)$ as a parameter that expresses the deformation of measurement associated with a measurement scale, $\Tr(f_y)$, we can create a sequence in which each successive deformation corresponds to a successive class of probability distributions with familial relations to each other defined by the structure of the sequence of deformations to $w(f_y)$.

\subsection{The general form of probability distributions}

From \Eq{shannonSoln}, the maximum entropy solution is
\begin{equation}\label{eq:shannonSoln2}
  p_y \propto m_ye^{-\Gl \Trf}.
\end{equation}
From \Eq{T_gen_soln}, we can now express the maximum entropy solution as
\begin{equation}\label{eq:master}
  p_y \propto m_y e^{-\GL e^{\Gb w}},
\end{equation}
where $\GL = \Gl(\Tr_0+\alpha/\Gb)$, and $w\equiv w(f_y)$. In the limit $\Gb \rightarrow 0$, this becomes
\begin{equation*}
  p_y \propto m_y e^{-\Gg w}
\end{equation*}
where $\Gg=\Gl\Ga$.

In Appendix B we describe the case of extreme values, for which we will use $m_y=\dd \Tr(f_y)/\dd y$.  When $f_y=y$ and $m_y=\dd\Tr(y)/\dd y=\Tr'$, it will be convenient to write
\begin{equation}\label{eq:Tprime}
  \Tr' \propto w'e^{\Gb w},
\end{equation}
where $w'=\dd w(y)/\dd y$, and as $\Gb\rightarrow0$, $\Tr' \propto w'$.

\section{Intuitive description of measurement and probability}

Intuitively, one can think of the symmetry of information invariance and measurement scale in the following way.  On a linear scale, each incremental change of fixed length yields the same amount of information or surprise independently of magnitude.  Thus, if we change the scale by multiplying all magnitudes by a constant, we obtain the same pattern of information relative to magnitude.  In other words, the linear scale is invariant to multiplication by a constant factor so that, within the framework of maximum entropy subject to constraint, we get the same information about probability distributions from an observation $y$ or $\GR(y)=cy$.  In this section, we use $f_y=y$ to isolate the symmetry expressed by particular choices of $\Tr$ and $\GR$.

On a logarithmic scale, each incremental change in proportion to the current magnitude yields the same amount of information or surprise.  Information is scale dependent.  We obtain the same information at any point on the scale by comparing ratios. For example, we gain the same information from the increment $\dd y/y=\dd\log(y)$ independently of the magnitude of $y$.  Thus, we achieve information invariance with respect to ratios by measuring increments on a logarithmic scale.  Within the framework of maximum entropy subject to constraint, we get the same information about probability distributions from an observation $y$ or $\GR(y)=y^c$, corresponding to informationally equivalent measurements $\Tr(y)=\log(y)$ and $\Tr(y^c)=c\log(y)$ (see ref.~\cite{frank10measurement}).

The form of a probability distribution under maximum entropy can be read directly as an expression of how the measurement scale changes with magnitude.  From the general solution in \Eq{shannonSoln}, linear scales $\Tr(y)\propto y$ yield distributions that are exponential in $y$, whereas logarithmic scales $\Tr(y)\propto c\log(y)$ yield distributions that are linear in $y^c$.  Exponential distributions of the form $e^{-\Gl y}$ arise from underlying linear scales, whereas power law distributions of the form $y^{-c}$ arise from underlying logarithmic scales.  

Many common distributions have compound form, in which one can read directly how the underlying measurement scale changes with magnitude.  For example, the gamma distribution has form $y^{-c}e^{-\Gl y}$.  When the magnitude of $y$ is small, the shape of the distribution is dominated by the power law component, $y^{-c}$.  As the magnitude of $y$ increases, the shape of the distribution is dominated by the exponential component, $e^{-\Gl y}$.  Thus, the underlying measurement scale grades from logarithmic at small magnitudes to linear at large magnitudes.  Indeed, the gamma distribution is exactly the expression of an underlying measurement scale that grades from logarithmic to linear as magnitude increases.  Nearly every common probability distribution can be read directly as a simple expression of the change in the underlying measurement scale with magnitude.

\section{Hierarchies of common probability distributions}

Given a particular form for the function $w(f_y)$, the measurement scale $\Tr(f_y)$ follows from \Eq{T_gen_soln} and the associated probability distribution follows from \Eq{master}.  Although we can choose $w$ in any way that we wish, certain measurement scales and information invariances are likely to be common.  We discussed in our earlier paper the importance two scales \cite{frank10measurement}.  The first scale grades from linear to logarithmic as magnitude increases, which we call the linear-log scale.  The second scale inverts the linear-log scale to be logarithmic at small magnitudes and linear at large magnitudes, giving the log-linear scale.  The inversion relating the two scales can be expressed by a Laplace transform, showing the natural duality of the scales and a connection to recent studies on superstatistics \cite{frank10measurement}. 

\subsection{The linear-log scale}

In terms of the notation in the present paper, we can define $w$ to establish a hierarchy of measurement deformations, in which each level in the hierarchy arises from successive application of the linear-log scaling to the scale in the previous level in the hierarchy.  

To define the linear-log measurement function in terms of $w$, note from \Eq{T_gen_soln} that, as $\Gb\rightarrow0$, the forms of $w$ and the measurement function $\Tr$ become equivalent with respect to setting the associated probability distribution. Thus, by setting $w$, we are defining the limiting form of the measurement function.  With these issues in mind, define
\begin{equation*}
  w^{(i)}=\log\left(c_i+w^{(i-1)}\right),
\end{equation*}
with $w^{(0)} = f_y$. The constant $c_i$ sets the transition between linear and logarithmic scaling: the scale is linear when $w^{(i-1)}$ is small relative to $c_i$ and logarithmic when $w^{(i-1)}$ is large relative to $c_i$.  As $c_i\rightarrow0$, we can use  $w^{(i)}=\log\left(w^{(i-1)}\right)$. 

It is easiest to see the abstract structure of the measurement hierarchy and the associated forms of probability distributions in the limiting case $c_i\rightarrow0$, leading to purely logarithmic deformations. The first row of Table~\ref{tab:pureLog} begins with the base measurement $w^{(0)}=f_y$. The following two rows show the first two deformations for the sequence $i=0,1,2$. 

\begin{table}[ht]
  \begin{center}
  \begin{tabular}{|l|l|l|}
    \hline
    $w(f_y)$ & 
    $p_y$ & 
    $p_{y|\beta \rightarrow 0} $ 
    \\
    \hline &&\\[-12pt]
    $f_y$ &
    $ m_y e^{-\GL e^{\Gb f_y}} $ & 
    $ m_y e^{-\Gg f_y} $ 
    \\
    $\log f_y$ & 
    $ m_y e^{-\GL f_y^\Gb} $ & 
    $ m_y f_y^{-\Gg} $
    \\
    $\log\log f_y$ & 
    $ m_y e^{-\GL\left(\log f_y\right)^\Gb} $ & 
    $ m_y\left(\log f_y\right)^{-\Gg} $ 
    \\
    \hline
  \end{tabular}
  \end{center}
  \caption{
  The logarithmic measurement hierarchy and the associated form of the probability distribution function $p_y$ from \Eq{master}. Note that $\Gb\rightarrow0$ of each line corresponds to $\Gb=1$ of the following line. 
  \label{tab:pureLog}
  }
\end{table}

This table gives the hierarchy of probability distributions that arise from successive logarithmic deformations.  With this structure in mind, we give the full expansion with $c_i\ne0$ in Table~\ref{tab:linearLog}.

\begin{table*}[ht]
  \begin{center}
  \begin{tabular}{|l|l|l|}
    \hline
    $w(f_y)$ & 
    $p_y$ & 
    $p_{y|\beta \rightarrow 0} $ 
    \\
    \hline &&\\[-12pt]
    $f_y$ &
    $ m_y e^{-\GL e^{\Gb f_y}} $ & 
    $ m_y e^{-\Gg f_y} $ 
    \\
    $\log \left(c_1+f_y\right)$ & 
    $ m_y e^{-\GL \left(c_1+f_y\right)^\Gb} $ & 
    $ m_y\left(c_1+f_y\right)^{-\Gg} $
    \\
    $\log\left(c_2+\log\left(c_1+ f_y\right)\right)$ & 
    $ m_y e^{-\GL\left(c_2+ \log\left(c_1+f_y\right)\right)^\Gb} $ & 
    $ m_y\left(c_2+ \log\left(c_1+f_y\right)\right)^{-\Gg} $ 
    \\
    \hline
  \end{tabular}
  \end{center}
  \caption{
  The linear-log measurement hierarchy.
  \label{tab:linearLog}
  }
\end{table*}

We discuss the interpretation of $m_y$ and $f_y$ below.  The different interpretations of those values lead directly to specific forms for probability distributions.  Before interpreting $m_y$ and $f_y$, we present an alternative measurement scale.  

\subsection{The log-linear scale}

We obtain the log-linear measurement deformation hierarchy \cite{frank10measurement} from 
\begin{equation*}
  w^{(i)}=c_i w^{(i-1)}+\log\left(w^{(i-1)}\right),
\end{equation*}
from which we obtain the probability distributions in Table~\ref{tab:logLinear}.  The log-linear scale changes logarithmically at small magnitudes and linearly at large magnitudes. 

\begin{table*}[ht]
  \begin{center}
  \begin{tabular}{|l|l|l|}
    \hline
    $w(f_y)$ & 
    $p_y$ & 
    $p_{y|\beta \rightarrow 0} $ 
    \\
    \hline &&\\[-12pt]
    $f_y$ &
    $ m_y e^{-\GL e^{\Gb f_y}} $ & 
    $ m_y e^{-\Gg f_y} $ 
    \\
    $c_1 f_y+\log f_y$ & 
    $ m_y e^{-\GL f_y^\Gb e^{c_1\Gb f_y}} $ & 
    $ m_y f_y^{-\Gg}e^{-c_1\Gg f_y} $
    \\
    $c_2\left(c_1f_y+\log f_y\right)+\log\left(c_1f_y+\log f_y\right)$ & 
    $ m_y e^{-\GL e^{\Gb w}} $ & 
    $ m_y e^{-\Gg w}$ 
%    $ m_y e^{-\GL\left(\left(c_1f_y+\log f_y\right)^\Gb+f_y^{c_2\Gb}+e^{c_1c_2\Gb f_y}\right)} $ & 
%    $ m_y\left(c_1f_y+\log f_y\right)^{-\Gg}f_y^{-c_2\Gg}e^{-c_1c_2\Gg f_y} $ 
    \\
    \hline
  \end{tabular}
  \end{center}
  \caption{
  The log-linear measurement hierarchy.  In the last line of the table, we use $w\equiv w(f_y)$ to shorten the expression.
  \label{tab:logLinear}
  }
\end{table*}

\subsection{Other scales}

The linear-log and log-linear scales describe common forms of measurement functions.  In this section, we briefly mention some other scales listed in Table~\ref{tab:otherScales}.  These additional scales illustrate the ways in which measurement relates to the patterns of probability.

The first line of Table~\ref{tab:otherScales} shows a log-linear-log scale for a measure on the interval $(c_1,c_2)$.  That scale changes logarithmically near the boundaries and linearly near the middle of the range, in which $\log b$ describes the skew in the scaling pattern.  

The second line of Table~\ref{tab:otherScales} shows a linear-log-linear scale for $f_y > 0$.  That scale changes linearly near the lower boundary of zero, linearly at large magnitudes, and logarithmically at intermediate values.  

\begin{table*}[ht]
  \begin{center}
  \begin{tabular}{|l|l|l|}
    \hline
    $w(f_y)$ & 
    $p_y$ & 
    $p_{y|\beta \rightarrow 0} $ 
    \\
    \hline &&\\[-12pt]
    $\log\left((c_2-f_y)(f_y-c_1)^b\right)$ &
    $ m_y e^{-\GL (c_2-f_y)^\Gb(f_y-c_1)^{b\Gb}} $ & 
    $ m_y (c_2-f_y)^{-\Gg}(f_y-c_1)^{-b\Gg} $ 
    \\
    $c_2f_y+b\log(c_1+f_y)$ &
    $ m_y e^{-\GL (c_1+f_y)^{b\Gb} e^{c_2\Gb f_y}} $ & 
    $ m_y (c_1+f_y)^{-b\Gg} e^{-c_2\Gg f_y} $ 
    \\
    \hline
  \end{tabular}
  \end{center}
  \caption{
  As $\Gb\rightarrow0$, line 1 is a log-linear-log measurement scale, and line 2 is a linear-log-linear measurement scale.
  \label{tab:otherScales}
  }
\end{table*}

\section{The scale of information}

The prior section presented probability distributions in terms of $m_y$ and $f_y$.  This section develops the interpretation of $m_y$, which arises from the relation between the scale of information invariance and the scale on which we express probability.

The key issue is that maximum entropy requires some underlying process to dissipate information.  With regard to deriving probability distributions, we may consider three aspects of scale in relation to the dissipation of information.  First, we may measure an outcome that arises from the aggregation of a series of random perturbations.  Second, we may measure only the extreme values of some underlying process, thereby throwing away all information about the underlying process except the form of the upper or lower tail of the underlying distribution.  Third, the dissipation of information may occur on one scale, but we may wish to make our measurements with respect to another scale.  

Each of these three aspects of the scale of information dissipation leads to a simple interpretation of probability measure in maximum entropy analysis.  We give a brief description each scale of information dissipation in relation to calculating $m_y$. 

\subsection{Aggregation of perturbations}

In the standard application of maximum entropy, accumulation of random perturbations without constraint leads to a uniform probability measure, which has maximum entropy and minimum information.  Thus, the scale at which information dissipates is the same as the scale of the probability measure.  In this case, our formulation of  maximum entropy has $m_y\equiv1$, because any information that arises from deformation of measurement relative to the uniform default is included in our expression of measurement scale, $\Tr(f_y)$.  

\subsection{Extreme values}

The distribution of extreme values depends only on the total (integral) of the probability measure in the tail of an underlying probability distribution \cite{embrechts97modeling}.  Because extreme value distributions arise from integrals of probability measures, the dissipation of information and the associated measurement scale for extreme values is expressed in terms of the cumulative distribution function (see Appendix B).  To obtain the associated form of the probability measure with respect to the probability distribution function, $p_y$, we must transform the invariant measurement scale originally expressed with respect to the integral of the underlying probability measure.  

To change from the integral scale of the cumulative distribution to the scale of the probability measure associated with the probability density function, we simply differentiate the initial measurement scale, $\Tr(f_y)$, from the cumulative distribution scale to obtain the associated change in probability measure (Appendix B).  For $f_y=y$, we obtain $m_y = \dd\Tr(y)/\dd y = \Tr'$.  We gave the general form of $m_y=\Tr'$ in \Eq{Tprime}.  

\subsection{Change of variable}

In some cases, information may dissipate on one scale, but we choose to express probability on another scale.  The log-normal distribution is the classic example.  Using Table~\ref{tab:pureLog}, we may consider measurements that lead to the Normal or Gaussian distribution by either analyzing squared deviations from a central value, $f_y=(y-\Gm)^2$ in line one of Table~\ref{tab:pureLog} with $\Gb\rightarrow0$ or, equivalently, linear perturbations of $f_y=(y-\Gm)$ in line two of Table~\ref{tab:pureLog} with $\Gb=2$.  In these cases, the perturbations are direct measures rather than the tail probabilities of extreme values, so $m_y=1$, and we have the standard form of the Gaussian as $p_y\propto e^{-\Gg(y-\Gm)^2}$.  

If we prefer to analyze values on a logarithmic scale, then we make the transformation $y\rightarrow\log y$.  This case does not arise from invariant information and the associated measurement transformation, but rather from a change of variable to a different scale.  So we must change our measure, as in any standard change of variable.  In this case, the change of measure is $m_y\dd y = \dd\log y=\dd y/y$, thus $m_y=y^{-1}$ and we obtain the log-normal distribution $p_y \propto y^{-1}e^{-\tilde{\Gg}\left(\log y - \tilde{\Gm} \right)^2}$, where $\tilde{\Gg}$ and $\tilde{\Gm}$ are transformed appropriately.

\section{Sufficiency: reduction of information}

The algorithm of maximum entropy allows us to choose any constraint $\Tr(f_y)$.  However, one of our main goals is to provide a clear rationale for the choice of constraint, so that maximum entropy is more than a tautological description of probability distributions.  We have expressed the choice of the measurement scale, $\Tr$, in terms of information invariance set by the underlying problem.  Although information invariance may take various forms, we followed our earlier paper \cite{frank10measurement} in which we defended the linear-log and log-linear scales as likely to be common scales associated with common information invariances.   

Once we have set the transformation $\Tr(f_y)$ by these common information invariances, many widely observed probability distributions follow.  In some cases, deriving probability distributions requires using an observable, $f_y\ne y$, that differs from the scale $y$ of the underlying probability measure.  For example, we may use the squared deviations from a central location, or a fractional moment $f_y=y^\Ga$, where $\Ga$ is not an integer \cite{frank09the-common}.  Use of $f_y=y$ or of squared deviations $f_y=(y-\Gm)^2$ is widely accepted.  Such choices lead to $f_y$ being a sufficient reduction of all of the information in observations in order to express common probability distributions.  

For our purposes in this paper, we simply note that we can derive many common distributions by the widely accepted use of $f_y=y$ or $f_y$ as a squared deviation.  But the reasons that particular choices of $f_y$ are sufficient have not been fully explained with regard to maximum entropy, particularly fractional moments such as $f_y=y^\Ga$ \cite{frank09the-common}. Those reasons probably have to do with the sort of analysis described by large deviation theory \cite{touchette09the-large}, in which the retained information arises from the minimal descriptions of location and scale that remain when one normalizes the consequences of a sequence of perturbations so that one obtains a stable limiting form.

 \begin{table*}[]
  \hbox{\null\hskip25pt
  {\small % reduce font size
  \begin{tabular}{|l|l|l|l|l|l|}
    \hline
    Distribution & 
    \hfil$p_y$\hfil &
    T.L.C & 
    $m_y$ &
    $f_y$ &
    Notes and alternative names
    \\
    \hline &&&&&\\[-12pt]
    Gumbel & 
    $e^{\Gb y - \GL e^{\Gb y}}$ &
    \ref{tab:pureLog}.1.2 & 
    $\Tr'$ &
    $y$ &
    \null
    \\
    Gibbs/Exponential & 
    $e^{-\Gg y}$ &
    \ref{tab:pureLog}.1.3 & 
    $\Tr',1$ &
    $y$ &
    \null
    \\
    Gauss/Normal & 
    $e^{-\Gg y^2}$ &
    \ref{tab:pureLog}.1.3 & 
    $1$ &
    $y^2$ &
    \null
    \\
    Log-Normal & 
    $y^{-1}e^{-\Gg \left(\log y\right)^2}$ &
    \ref{tab:pureLog}.1.3 & 
    $y^{-1}$ &
    $y^2$ &
    Change of variable $y\rightarrow\log y$
    \\
    Fr\'echet/Weibull& 
    $y^{\Gb-1}e^{-\GL y^\Gb}$ &
    \ref{tab:pureLog}.2.2 & 
    $\Tr'$ &
    $y$ &
    \null
    \\
    Stretched exponential & 
    $e^{-\GL y^\Gb}$ &
    \ref{tab:pureLog}.2.2 & 
    $1$ &
    $y$ &
    Gauss with $\Gb = 2$
    \\
    Symmetric \levy & 
    $e^{-\GL |y|^\Gb}$ (Fourier domain) &
    \ref{tab:pureLog}.2.2 & 
    $1$ &
    $|y|$ &
    $\Gb \le 2$; Gauss ($\Gb=2$), Cauchy ($\Gb=1$); ref.~\cite{frank09the-common}
    \\
    Pareto type I& 
    $y^{-\Gg}$ &
    \ref{tab:pureLog}.2.3 & 
    $\Tr',1$ &
    $y$ &
    \null
    \\
    Log-Fr\'echet& 
    $y^{-1}(\log y)^{\Gb-1}e^{-\GL(\log y)^\Gb}$ &
    \ref{tab:pureLog}.3.2 & 
    $\Tr'$ &
    $y$ &
    Also from Fr\'echet: $y\rightarrow\log y$, $m_y=y^{-1}\Tr'(y)$
    \\
    ?? & 
    $e^{-\GL(\log y)^\Gb}$ &
    \ref{tab:pureLog}.3.2 & 
    $1$ &
    $y$ &
    Also stretched exponential with $f_y=\log y$
    \\
    Log-Pareto type I & 
    $y^{-1}\left(\log y\right)^{-\Gg-1}$ &
    \ref{tab:pureLog}.3.3 & 
    $\Tr'$ &
    $y$ &
    Log-gamma; Pareto I: $y\rightarrow\log y$, $m_y=y^{-1}$
    \\
    ?? & 
    $\left(\log y\right)^{-\Gg}$ &
    \ref{tab:pureLog}.3.3 & 
    $1$ &
    $y$ &
    Also from Pareto I with $f_y=\log y$
    \\
    Pareto type II & 
    $\left(c_1+y\right)^{-\Gg}$ &
    \ref{tab:linearLog}.2.3 & 
    $1$ &
    $y$ &
    Lomax
    \\
    Generalized Student's & 
    $\left(c_1+y^2\right)^{-\Gg}$ &
    \ref{tab:linearLog}.2.3 & 
    $1$ &
    $y^2$ &
    Pearson type VII, Kappa; includes Cauchy
    \\
    ?? & 
    $\left(\log\left(c_1+y\right)\right)^{-\Gg}$ &
    \ref{tab:linearLog}.3.3 & 
    $1$ &
    $y$ &
    $c_2=0$; also Pareto I with $f_y=\log (c_1+y)$
    \\
    Gamma & 
    $y^{-\Gg}e^{-c_1\Gg y}$ &
    \ref{tab:logLinear}.2.3 & 
    $1$ &
    $y$ &
    Pearson type III, includes chi-square
    \\
    Generalized gamma & 
    $y^{-k\Gg}e^{-c_1\Gg y^k}$ &
    \ref{tab:logLinear}.2.3 & 
    $1$ &
    $y^k$ &
    Chi with $k=2$ and $c_1\Gg=1/2$
    \\
    Beta & 
    $(c_2-y)^{-\Gg}(y-c_1)^{-b\Gg}$ &
    \ref{tab:otherScales}.1.3 & 
    $1$ &
    $y$ &
    Pearson type I; log-linear-log on $(c_1,c_2)$
    \\
    Beta prime/F& 
    $y^{-b\Gg}(1+y)^{(b+1)\Gg}$ &
    \ref{tab:otherScales}.1.3 & 
    $1$ &
    $\frac{y}{1+y}$ &
    Pearson type VI, $y>0$
    \\
    Gamma variant& 
    $(c_1+y)^{-b\Gg} e^{-c_2\Gg y}$ &
    \ref{tab:otherScales}.2.3 & 
    $1$ &
    $y$ &
    Linear-log-linear pattern as $y$ rises from zero
    \\
   \hline
  \end{tabular}
  } % end small font size
  } % end hbox
  \caption{
  Some common probability distributions.  The column T.L.C gives the table, line, and column of the underlying form presented in the earlier tables of abstract distributions. For example, \ref{tab:pureLog}.1.2 refers to Table~\ref{tab:pureLog}, first line, second column.  The measurement adjustment is given as either $m_y=1$ for direct scales, or $m_y=\Tr'$ for extreme values as in \Eq{Tprime}, along with any consequences from a change of variable such as $y\rightarrow\log y$.  Cases in which the same structural form arises for either $m_y=\Tr'$ or $m_y=1$ are shown as $\Tr',1$, without adjusting parameters for trivial differences.  The value of $f_y$ gives the reduction of data to sufficient summary form. Direct values $y$, possibly corrected by displacement from a central location, $y-\Gm$, are shown here as $y$ without correction.  Squared deviations $(y-\Gm)^2$ from a central location are shown here as $y^2$. See refs.~\cite{johnson94continuous,johnson95continuous,kleiber03statistical} for listing of distributions.  Many additional forms can be generated by varying the measurement function.  In the first column, the question marks denote a distribution for which we did not find a commonly used name.
  \label{tab:commonDistn}
  }
\end{table*}

\section{Conclusions}

Table~\ref{tab:commonDistn} shows many of the commonly observed probability distributions.  Those distributions arise directly from maximum entropy applied to various natural measurement scales.  The measurement scales express information invariances associated with particular types of problems and the scale on which information dissipation occurs.  We confined ourselves to various combinations of linear and logarithmic scaling, which were sufficient to express many common distributions.  Our method readily extends to other types of information invariance and measurement scale and their associated probability distributions.  

\section*{Acknowledgements}

SAF is supported by National Science Foundation grant EF-0822399,
National Institute of General Medical Sciences MIDAS Program grant
U01-GM-76499, and a grant from the James S.~McDonnell Foundation.
DES thanks Insight Venture Partners for support.

\newpage

\section*{Appendices}

\appendix

\section{On the association between measurement functions and classes
  of scale transformations}

If the transformation $f_y \rightarrow G \! \left( f_y \right)$ is an invariance of a measurement function $T$, it is clear that repeated applications of $G$, expressed as $G \circ G, G \circ G \circ G,\ldots$, are also invariances of $T$.  It is the larger group of invariances that we wish to identify with the measurement scale that defines $T$, and not only a single transformation.  To simplify notation in this Appendix, we use $f_y=y$.  The same analysis applies to $f_y$.

In general, making a unique association between a transformation $G$ and a measurement function $T$ is inconvenient for finite transformations, because $G$ combines a magnitude and a direction of deformation.  The magnitude is added under compositions $G \circ G \ldots$, while the direction remains invariant.  As we will derive below, the relevant measure of the magnitude of a transformation as in \Eq{transDefInv} will be $\sim \log b$, and the relevant measure of direction will be $a / \left( b-1 \right)$.  To isolate the direction of $G$ that may be associated with a measurement function $T$, we work with infinitesimal rather than finite affine transformations.

Infinitesimal transformations are constructed from \Eq{transDefInv} in the text by writing $a \equiv \epsilon \alpha$, $\left( b-1 \right) \equiv \epsilon \beta$, and then taking $\epsilon \rightarrow 0$ for fixed $\alpha$ and $\beta$.  An infinitesimal transformation $G^{\epsilon}$ then satisfies \Eq{transDefInv} in the form
\begin{equation}
  T \!
  \left[ 
    G^{\epsilon} \!
    \left( y \right)
  \right] = 
  T \! \left( y \right) + 
  \epsilon 
  \left[ 
    \alpha + 
    \beta 
    T \! \left( y \right)
  \right] . 
\label{eq:T_G_inf_def}
\end{equation}
$G$ itself must therefore also be infinitesimally different from the identity, and must have the form
\begin{equation}
  G^{\epsilon} \!
  \left( y \right) = 
  y + 
  \epsilon 
  v \! \left( y \right) . 
\label{eq:G_inf_on_y}
\end{equation}
for some function $v \! \left( y \right)$.

We introduce a quantity $\hat{v}$ called the \emph{generator} of the deformation, such that the operator $e^{\epsilon \hat{v}}$ generates the infinitesimal transformations Eqs.~(\ref{eq:T_G_inf_def},\ref{eq:G_inf_on_y}), and such that finite transformations $G$ or affine transformations \Eq{transDefInv} are produced by the exponential operation of $\hat{v}$ with non-infinitesimal $\epsilon$.  Compounding a function corresponds to addition of parameters $\epsilon$, as may be checked from the power-series definition of $e^{\epsilon \hat{v}}$ within its radius of convergence.

We define a \emph{representation} of the generator $\hat{v}$ as an explicit differential operator that produces the correct transformation on the argument $y$ or $T \! \left( y \right)$, as appropriate.  The two representations of the generators are related as
\begin{eqnarray}
  T \!
  \left[ 
    y + \epsilon 
    v \! \left( y \right)
  \right] 
& = & 
  \left[ 
    1 + \epsilon 
    v \! \left( y \right)
    \frac{\dd}{\dd y}
  \right] 
T \! \left( y \right)
\nonumber \\
& = & 
  \left[ 
    1 + \epsilon 
    \left( 
      \alpha + 
      \beta T 
    \right)
    \frac{\dd}{\dd T}
  \right] 
  T . 
\nonumber \\
\label{eq:T_G_inf_reps}
\end{eqnarray}
From the requirement that the two expressions produce the same result, we may assign the representations
\begin{eqnarray}
  \hat{v} 
& \leftrightarrow & 
  v \! \left( y \right)
  \frac{\dd}{\dd y} \equiv 
  \frac{\dd}{\dd w}
\nonumber \\
& \leftrightarrow & 
  \left( 
    \alpha + 
    \beta T 
  \right)
  \frac{\dd}{\dd T}  
\label{eq:v_reps}
\end{eqnarray}
for some function $w \! \left( y \right) = \int^y \dd y^{\prime} 1 / v \!
\left( y^{\prime} \right)$. 

Regarding $T$ as a function of argument $w$ rather than $y$, and setting equal the two coefficients of $\epsilon$ in \Eq{T_G_inf_reps}, we obtain a relation between any function $w \! \left( y \right)$, coefficients $\alpha$ and $\beta$, and the function $T$ in the form
\begin{equation}
  \frac{\dd T}{\dd w} = 
  \alpha + \beta 
  T . 
\label{eq:reps_w_to_diff_eq_T}
\end{equation}
This is rearranged to produce \Eq{T_resp_genform}.

From the solutions to \Eq{reps_w_to_diff_eq_T}, we may readily check that the action of the transformation $e^{\epsilon \hat{v}}$ for arbitrary $\epsilon$ (not necessarily small) is
\begin{equation}
  T \! 
  \left[ 
    G \! 
    \left( y \right)
  \right] = 
  e^{\epsilon \hat{v}}
  T \! \left( y \right) = 
  \frac{\alpha}{\beta}
  \left( 
    e^{\epsilon \beta} - 1 
  \right) + 
  e^{\epsilon \beta} 
  T \! \left( y \right) , 
\label{eq:finite_trans_T}
\end{equation}
from which we recover expressions for the coefficients $a$ and $b$ in \Eq{transDefInv}.  Under composition $G \rightarrow G \circ G$, the parameter $\epsilon \rightarrow 2 \epsilon$.  The composition rules for $a$ and $b$ under composition of $G$ may be worked out easily, but depending on the function $w \! \left( y \right)$, the direct composition of finite transformations $G$ on $y$ may be quite complicated.

\section{Information measures for cumulative distributions}
\label{sec:cum_dist}

The presence of the measure $m_y$ in the probability density function in \Eq{shannonSoln2} complicates the discussion of measurement invariance, because in the general case $m_y$ is not required to obey any prescribed transformation when $f_y \rightarrow G \! \left( f_y \right)$.  In general, $y$ need not even be a numerical index, whereas $T \! \left( f_y \right)$ is necessarily numerical because it is proportional to an information measure $- \log \left( p_y / m_y \right)$.

The class of cases in which the measurement function, $T$, completely controls the properties of $p_y$ are those in which measurement constrains the \emph{cumulative} probability distribution function rather than the probability density function.  For these cases $m_y$ is not independent, but is given in terms of $T$ and $f_y$, as we now show.

Relative entropy is ordinarily defined for the probability density.  However, if we set
\begin{equation}
  m_y = 
  \frac{
    \dd T \! \left( f_y \right)
  }{
    \dd y 
  } = 
  \frac{\dd f_y}{\dd y}
  T^{\prime} \! 
  \left( f_y \right) , 
\label{eq:m_y_cum_ent}
\end{equation}
then $m_y$ becomes a Lebesgue measure on $y$ with respect to the increment $\dd T$.  The probability density from \Eq{shannonSoln2} becomes
\begin{equation}
  p_y \propto
  \frac{\dd }{\dd y}
  e^{- 
    \lambda T  \left( f_y \right)
  } . 
\label{eq:dens_from_cum}
\end{equation}
\Eq{dens_from_cum} defines the relation between a probability density and its cumulative distribution, meaning that under a suitable ordering of $y$, we may take $e^{- \lambda T\left( f_y \right) }$ to be the cumulative distribution.

With this choice of measure, the relative entropy $\mathcal{E}$ from \Eq{shannonDef} becomes
\begin{widetext}
\begin{eqnarray}
  - \int \dd y\, p_y 
  \log 
  \left( \frac{p_y}{m_y} \right) 
& = & 
  - \int \dd y 
  \frac{\dd T}{\dd y}
  \left( \frac{p_y}{\dd T / \dd y} \right) 
  \log 
  \left( \frac{p_y}{\dd T / \dd y} \right)   
\nonumber \\
& = & 
  - \int \dd T
  p_T 
  \log 
  p_T , 
\label{eq:rel_ent_T}
\end{eqnarray}
\end{widetext}
in which $p_T$ is the probability density defined on the variable $T$. Since the maximum-entropy solution is always exponential in $T$, the relative entropy of \Eq{rel_ent_T} is effectively an information function for the cumulative distribution.

An application in which constraints under aggregation apply by construction to the cumulative distribution is the computation of extreme-value statistics~\cite{kotz00extreme}.  The cumulative probability distribution for the maximizer or minimizer of a sample of $n$ realizations of a random variable is the product of $n$ factors of the cumulative distribution for a single realization.

It was also noted in ref.~\cite{frank09the-common} that the relative entropy may be evaluated on the characteristic function (Fourier or Laplace transform) of a distribution, and that the maximum-entropy solutions in the transformed domain are the L\'{e}vy stable distributions.  The characteristic function at frequency argument $k = 0$ always takes value unity.  Therefore it, like a cumulative distribution, has a reference normalization of unity, and indeed, the symmetric L\'{e}vy-stable distributions~\cite{sato01basic} correspond in form to the Weibull family of extreme value distributions.  Both are obtained within our classification for $m_y$ defined by \Eq{m_y_cum_ent}, for suitable reductions $f_y$.

\section{The Morris Natural Exponential Families in relation to entropy-maximizing distributions}

\subsection{Symmetry-based approaches to select or to classify probability distributions} 

Many systems, since Pearson's, for either selecting or classifying
probability distributions, have been based on symmetry groups, as our
method is.  (Pearson's system may be seen as one based on the analytic
structure of the log-probability, a criterion that we will return to
consider in a moment.)  The systems differ in generality, depending on the space in which the symmetry group acts, and depending on whether it constrains a single distribution or a family.  
Two methods based on symmetry (ours and that of
Carl Morris, described below) have interpretations in terms of 
\emph{scale invariance} of observables.  Both systems collect probability distributions into families, whose members differ only by a scale factor.  A third approach (known as \emph{Objective Bayesian} methods) applies symmetry to
the underlying measure space which, as we note in Appendix~\ref{sec:cum_dist}, may be very different from the space of observed magnitudes.  This approach is concerned not directly with
families of distributions, but with the particular distribution
defined by a reference measure.  We will briefly summarize the overlaps and differences of these methods.  

Objective Bayesian methods, initiated by Jeffries~\cite{jeffries57scientific} but given the interpretation of objectivity largely by Jaynes~\cite{jaynes68prior}, recognize that the reference measure $m_y$ in a relative entropy---beyond being needed to make logarithms well-defined and independent of change of variables---may reflect information about measurement scales.  By ensuring that the reference measure is consistent with known symmetries of the phenomenon under study (which are not generally expressed within particular sample observations), Objective Bayesian methods seek to systematize the entire maximum-entropy procedure.  This use of the reference measure is consistent with our treatment of measurement, though by itself it is more limited, as we discuss in Ref.~\cite{frank10measurement}, and it may also be misleading in cases~\cite{seidenfeld79protect}.  In the context of the present discussion, the most important limitation of Objective Bayesian methods is that they select properties of a single distribution $m_y$, rather than properties of a family.  

Our approach broadens the class of symmetries that can be 
considered, beyond those available to Objective Bayesian methods, as discussed in Sec.~7 of Ref.~\cite{frank10measurement}.  Through the measurement function, it relates a potentially nonlinear contour of deformations of measured magnitudes to a linear transformation within the affine group that exists for general maximum-entropy problems.  We have
embedded distributions within a hierarchy by using the two-parameter
freedom of the affine group to provide a range
of responses of information to the change in the scale of measurement.

\subsection{The Morris classification of distributions in relation to
  maximum entropy}

In a pair of papers in 1982 and
1983~\cite{morris82natural,morris83natural}, Carl Morris
proposed another classification system for probability distributions,
which overlaps both with Objective Bayesian methods and with our
approach.  Like our method, Morris's concerns families of
probability distributions generated by a change in constraint or
measurement scale.  Like all of the approaches we have mentioned,
Morris's system uses relative entropy in a conventional maximization framework.  That system differs from ours in using only a linear constraint on what Morris terms the \emph{natural observation}, and obtaining nonlinear dependence on that constraint through a second boundary condition placed on entropy.  

The Morris system blends interesting elements of Pearson's
restrictions on analytic structure, our use of symmetry, and the
Objective Bayesian concern with the reference measure, as follows:
Morris considers distribution families that are invariant under offset and rescaling of the natural observation, which Morris labels $X$, and which is analogous to using a coordinate system that is always linear in our $f_y$.   His classification therefore does not  not invoke any explicit representation of the symmetries inherent in differing measurement systems.  In order to encompass distributions that are not simply exponential in the values $x$ (taken by the observation $X$), he instead restricts the form of the reference measure in a relative entropy, analogous to our $m_y$.  Unlike Objective Bayesian methods, however, this restriction does not come from the direct action of a symmetry on the reference measure, but rather from the form of the relative entropy across the family of distributions produced by scale change.  

The classification system of Ref.~\cite{morris82natural} derives from the
cumulant-generating function and the relation between the variance
and the mean as the parameter in this generating function is shifted.
The distributions that define the cumulant-generating function
constitute what Morris calls \emph{natural exponential families} (NEF),
and the dependence of variance on mean within these families is
restricted in his system to be an exact quadratic polynomial.  The resulting subclass of distributions within the
NEF class is termed QVF (for \emph{quadratic variance function}).  The mean-variance relation that defines the NEF-QVF distributions is preserved
under offset and rescaling of the natural observation, and under convolution.  Therefore, the distributions in this class would be
expected to arise frequently in problems of aggregation.  We show in this appendix that the QVF condition is equivalent to the
requirement that the relative entropy over a family has the form (up
to analytic continuation) of a Kullback-Leibler divergence.  The analytic continuation is 
determined by the roots of the quadratic variance polynomial, and these roots in turn have a relation to the roots for log-probability in the Pearson system. 

The distributions selected by Morris's criterion are either bounded,
or have exponential or faster decay in their tails.  We show that, when they are classified according to
their analytic structure, they are in fact either interior members or
degenerate limits of only two families of distributions: One family of
continuous-valued distributions is associated with complex-conjugate roots of the variance function, and a complex 
analytic continuation of the Kullback-Leibler form for relative
entropy.  A second family of discrete-valued distributions is associated with
real-valued roots, and real-valued continuations of the Kullback-Leibler relative entropy.  In this sense, the Morris classification
shows that six important distribution families are in fact selected by
a single set of invariances---of these, the offset and scale invariances are instances of our linear measurement rescaling.   These selected families are therefore very commonly observed,
but also rather tightly restricted.  Preservation of a functional class under convolution is similar to the criterion leading to the extreme-value or L\'{e}vy distributions, as we have discussed in the main text, and is therefore one of many  forms of measurement invariance that may be considered.  

Here we will re-formulate the Morris criterion and its solutions within a standard framework of maximum entropy.  We will show that the role of the reference measure in a relative entropy is equivalent to that of a \emph{second} observed quantity, which will generally be linearly independent of the natural observation $X$.  Scale change of the natural observation defines what is known as an expansion path, which consists of the distributions within an exponential family.  The second observed quantity, associated with the reference measure, is given a gradient constraint rather than a value constraint.  It is through the interaction of these two constraints that nonlinear dependence on $x$ is obtained in the log-probability.  At the end of the Appendix we mention a relation between the Morris system and the Pearson system based on the log-probability.   When the Morris QVF criterion is expressed as a formal constraint on
entropy, this form is imposed on the leading terms of log-probability
by the large-deviations property of cumulant-generating functions.

\subsection{Definition of the natural exponential families}

The NEF distributions are defined in relation to the cumulant-generating function, which arises naturally in the method of maximum entropy.  The most direct way to re-formulate the original presentation of Refs.~\cite{morris82natural,morris83natural} in terms of maximum entropy is to assume a (Shannon-type) entropy in a higher-dimensional state space than the univariate space of the natural observation $X$.  The high-dimensional states have non-uniform density when they are projected onto the one dimension in which the probability distribution varies.  Once a Lagrangian is defined from this initial re-formulation, it becomes easy to re-interpret the density of states as a reference measure in a relative entropy (and the latter interpretation is more general).  The cumulant-generating function is then the Legendre transform of
this relative entropy.  We develop the two interpretations in order,
to connect the derivations of Refs.~\cite{morris82natural,morris83natural}
systematically to the formulation we use in the main text.

\subsubsection{The Stieltjes measure as a density of states}

Ref.~\cite{morris82natural} introduces a Stieltjes measure $\dd F \! \left( x
\right)$, and an initial probability distribution $P_0$ associated with this measure, defined by 
\begin{equation}
  P_0 \! 
  \left(
    X \in A 
  \right) = 
  \int_A 
  \dd F \! \left( x \right) , 
\label{eq:Morris_pdf_0}
\end{equation}
for an arbitrary set $A$ in the range of $x$.  
With respect to this original probability measure, Morris introduces the
exponential families in terms of a probability mass function 
\begin{equation}
  \phi \! \left( x \mid \theta \right) = 
  e^{
    x \theta - 
    \psi \left( \theta \right)
  } , 
\label{eq:Morris_solution}
\end{equation}
which multiplicatively weights the original measure $\dd F \! \left( x
\right)$.  The normalizing constant $\psi \! \left( \theta \right)$ in
Eq.~(\ref{eq:Morris_solution}) is the cumulant-generating function,
given by
\begin{equation}
  e^{
    \psi \left( \theta \right)
  } \equiv 
  \int 
  \dd F \! \left( x \right)
  e^{x \theta} .  
\label{eq:Morris_CGF}
\end{equation}
The NEF distributions are the normalized versions of
the distributions that define the cumulant-generating function.
In the original Stieltjes measure, the probabilities defined from these distributions are 
\begin{equation}
  P \! 
  \left(
    X \in A 
  \right) = 
  \int_A 
  \dd F \! \left( x \right)
  e^{
    x \theta - 
    \psi \left( \theta \right)
  } .  
\label{eq:Morris_pdf}
\end{equation}

With respect to the measure $\dd F \! \left( x \right)$, we may
obtain the solutions~(\ref{eq:Morris_solution}) by extremizing the
Lagrangian 
\begin{widetext}
\begin{equation}
  \mathcal{L} = 
  - \int
  \dd F \! \left( x \right) \, 
  \phi \! \left( x \right) \log \phi \! \left( x \right) + 
  \theta
  \left(
  \int
    \dd F \! \left( x \right) \, 
    \phi \! \left( x \right) x - 
    \mu
  \right) - 
  \kappa 
  \left( 
  \int
    \dd F \! \left( x \right) \, 
    \phi \! \left( x \right) - 
    1
  \right) 
\label{eq:Lagrange_Morris}
\end{equation}
\end{widetext}
over its natural argument $\phi \! \left( x \right)$ and the Lagrange
multipliers $\theta$ and $\kappa$.  Here we have replaced the 
notation $\lambda$ from the text with Morris's $\theta$ for ease of
reference.  From its role as a normalization constant, the multiplier
$\kappa$ must evaluate to the cumulant-generating function $\psi \!
\left( \theta \right)$ on solutions.

Lagrangian problems of this form arise frequently in systems where a
high-dimensional state space is projected down onto a single
coordinate $x$, which is the only observed property on which
distributions depend.  The Lagrangian~(\ref{eq:Lagrange_Morris})
effectively treats $\phi \! \left( x \right)$ as the ratio of a
probability density to a \emph{uniform} reference measure on the
original high-dimensional space.  The Stieltjes measure $\dd F \! \left(
x \right)$ is the marginal projection of the original measure onto the
coordinate $x$, and the derivative $\dd F / \dd x$ is known as the
\emph{density of states}.  ($\dd F \! \left( x \right)$ need not be
smooth, and $\dd F / \dd x$ may readily be a non-continuous distribution,
such as a sum of Dirac $\delta$-functions, representing a discrete
rather than continuous probability density).

The entropy in this formulation appears as a standard Shannon entropy
(equivalent to a relative entropy with a uniform reference measure) in the
high-dimensional coordinates.  It evaluates to the Legendre transform
of the cumulant-generating function, 
\begin{eqnarray}
  S \! 
  \left( 
    \mu \! \left( \theta \right)
  \right) 
& \equiv &  
      - \int
      \dd F \! \left( x \right) \, 
      p \! \left( x  \mid \theta \right) 
      \log p \! \left( x  \mid \theta \right) 
\nonumber \\
& = & 
  - \theta \mu \! \left( \theta \right) + 
  \psi \! \left( \theta \right) , 
\label{eq:Morris_ent_eval}
\end{eqnarray}
in which $\mu  \! \left( \theta \right)$ is the mean value in the distribution $p \! \left( x  \mid \theta \right)$.   $\theta$ is the natural argument of $\psi$, while $\mu$ from the
variational problem is the natural argument of $S$.  Therefore it is
usual to write this Legendre transform pair as
\begin{eqnarray}
  \psi \! \left( \theta \right)
& = & 
  {
    \left. 
      S \! \left( \mu \right) - 
      \mu
      \frac{\dd S}{\dd\mu}
    \right| 
  }_{
    \mu = \mu \left( \theta \right)
  }
\nonumber \\
  S \! \left( \mu \right)
& = & 
  {
    \left.
      \psi \! \left( \theta \right) - 
      \theta
      \frac{\dd\psi}{\dd\theta}
    \right|
  }_{
    \theta = \theta \left( \mu \right)
  }
\label{eq:Lengendre_pair}
\end{eqnarray}
In the second line, $\theta \! \left( \mu \right)$ is the inverse
function to $\mu \! \left( \theta \right)$.  (In statistical
mechanics, where $-\theta$ is the inverse temperature if $x$ is the
energy, $\psi$ arises as $\theta$ times the \emph{Helmholtz Free
  Energy}.)

We note several properties of these functions that will be useful in
understanding Morris's NEF-QVF families.  When $\theta = 0$ no
correction to the normalization is needed in $P \!  \left( X \in A
\right)$, so we have immediately that $\psi \! \left( 0 \right) = 0$ as
well.  If we denote by ${\mu}_0 \equiv \mu \! \left( 0 \right)$, then
it follows that $S \! \left( {\mu}_0 \right) = 0$ also.  The
definition of the Legendre transform pair~(\ref{eq:Lengendre_pair})
gives the important dual relations 
\begin{eqnarray}
  \frac{
    \dd \psi \! \left( \theta \right)
  }{
    \dd \theta 
  }
& = & 
  \mu \! \left( \theta \right)
\nonumber \\
  \frac{
    \dd S \! \left( \mu \right)
  }{
    \dd \mu 
  }
& = & 
  - \theta \! \left( \mu \right) . 
\label{eq:Lengendre_pair_grads}
\end{eqnarray}
It follows that ${\left. \dd S / \dd\mu \right|}_{{\mu}_0} = 0$.  With
these two constants of integration, $S \! \left( \mu \right)$ will be
completely specified by the form of its second derivative.

\medskip
\subsubsection{Replacing the density of states with a reference measure in relative entropy}

For the univariate distributions, whether continuous or discrete, we
may define a shorthand for Eq.~(\ref{eq:Morris_pdf}) by identifying
the probability density function on $x$ as
\begin{equation}
  p_{x \mid \theta} \equiv 
  \frac{\dd F}{\dd x}
  e^{
    x \theta - 
    \psi \left( \theta \right)
  } .  
\label{eq:p_x_ident}
\end{equation}
The Lagrangian~(\ref{eq:Lagrange_Morris}) becomes, under this change
of variable, 
\begin{widetext}
\begin{equation}
  \mathcal{L} = 
  - \int
  \dd x \, 
  p_x \log 
  \left(
    \frac{p_x}{\dd F/\dd x}
  \right) + 
  \theta
  \left(
  \int
    \dd x \, 
    p_x x - 
    \mu
  \right) - 
  \kappa 
  \left( 
  \int
    \dd x \, 
    p_x - 
    1
  \right) . 
\label{eq:Lagrange_Morris_KL}
\end{equation}
\end{widetext}
The constraint terms are unchanged, but the entropy is now manifestly
a \emph{relative entropy} for the density $p_x$ with reference measure
$\dd F/\dd x$.  

\subsubsection{Arriving at nonlinear expansion paths through mixed
  boundary conditions} 

The Morris families, like the Pearson families and like our classes
based on measurement, include distributions that are nonlinear in the
values $x$ taken by the natural observation $X$.  Both Morris's
families and ours are based on affine transformation, so that their
distributions form what are known as \emph{expansion paths}.  (This
term is used also in economics for constrained maximization problems,
in which $\mu$ generally describes a budget constraint.  The original
usage, in statistics, is mentioned in
Ref.~\cite{morris82natural}.)  Whereas we achieve nonlinear
dependence on $x$ by considering the symmetries of measurement, the
Morris system achieves nonlinearity through the use of mixed boundary
conditions, when this system is described in terms of entropy
maximization.  By using two constraints---one to specify the family
and the other to fix a point on the expansion path---Morris is able
to apply a fixed-gradient condition with respect to one constraint,
and a fixed-value condition for the natural observation.  Because we
specify distributions from the affine transformation of a single
observable, we must incorporate nonlinearities into the measurement
function itself.

Here, we derive the NEF criterion by converting the relative
entropy to a form in which the reference measure may be interpreted as
a second observable.  The ubiquitous use, in statistical physics and
thermodynamics, of cumulant-generating functions and their Legendre
transforms under mixed boundary conditions, provides intuition from
familiar systems for the meaning of the resulting expansion paths.  In
the next section we derive the way in which the QVF condition of
Morris then places constraints on the reference measure, which plays
the role of the secondary observation. 

The Lagrangian~(\ref{eq:Lagrange_Morris_KL}) is an instance of a more
general class of maximum entropy problems in which the relative entropy has
uniform measure (and therefore has the form of a Shannon entropy), and
the reference measure appears as an additional constraint term,
\begin{widetext}
\begin{equation}
  \mathcal{L} = 
  - \int
  \dd x \, 
  p_x \log p_x + 
  \theta
  \left(
  \int
    \dd x \, 
    p_x x - 
    \mu
  \right) + 
  \lambda
  \int
    \dd x \, 
    p_x 
    \log 
    \left( 
      \frac{\dd F}{\dd x}
    \right) - 
  \kappa 
  \left( 
  \int
    \dd x \, 
    p_x - 
    1
  \right) . 
\label{eq:Lagrange_Morris_KL_grad}
\end{equation}
\end{widetext}
Here a variable $\lambda$ has been added as a \emph{parameter} in the
variational problem, parallel to the parameter $\mu$ in the constraint
on $\int \dd x \, p_x x$.  When $\lambda = 1$,
Eq.~(\ref{eq:Lagrange_Morris_KL_grad}) reduces to
Eq.~(\ref{eq:Lagrange_Morris_KL}).  and the choice of reference
measure does not matter because it cancels in the two logarithms.  For
more general $\lambda$, a uniform reference measure is explicitly
required to make the logarithms well-defined.  The distribution
solving Eq.~(\ref{eq:Lagrange_Morris_KL_grad}) is
\begin{equation}
  p_x = 
  e^{
    \lambda 
    \log \left( \text{d} F/\text{d} x \right) + 
    \theta x - 
    \kappa
  } . 
\label{eq:solving_gen_lamb}
\end{equation}

The Shannon entropy of Eq.~(\ref{eq:Lagrange_Morris_KL_grad}) is
maximized subject to mixed constraints, which may be seen as follows.
The entropy with two constraint terms is a function of two arguments
$S \! \left( \mu, \xi \right)$, where $\xi = \left< \log \left( \dd F/\dd x
\right) \right>$ at the given values of $\lambda$ and $\mu$.  Then $\lambda = -
\partial S / \partial \xi$, just as $\theta = - \partial S / \partial
\mu$ from Eq.~(\ref{eq:Lengendre_pair_grads}).  Because $\mu$ is an
argument to the entropy, whereas $\lambda$ is a gradient, problems of
this sort resemble solutions to differential equations under mixed
Dirichlet and Neumann boundary conditions. 

The set of distributions~(\ref{eq:solving_gen_lamb}), as $\lambda$ is
held fixed and $\mu$ is varied, make up the expansion path for the
entropy with respect to constraint $\int \dd x \, p_x x$.  The natural
exponential families are the distributions on this expansion path,
given a gradient constraint with respect to the observable $\int \dd x \,
p_x \log \left( \dd F / \dd x \right)$.

\subsection{The subset of natural exponential families with quadratic
  variation} 

Any reference measure may in principle form the basis for an expansion
path with mixed constraints.  In contrast to Objective Bayesian
methods, in which $\log \left( \dd F / \dd x \right)$ is constrained by
symmetry, the Morris system constrains reference measures by
restricting the form of the variance function---equivalent to
restricting the form of the \emph{entropy}---along the nonlinear
expansion path. 

\subsubsection{The QVF family and Kullback-Leibler entropies}

The definition of the cumulant-generating function is that, not only
does $\dd\psi / \dd\theta = \mu$, but $\dd^2 \psi / {\dd \theta}^2$ is the
variance of the observation $X$.  Morris defines its relation to the
mean $\mu$ as a \emph{variance function} $V \!  \left( \mu \right)$.
The \emph{quadratic variance relation} is the dependence
\begin{equation}
  \frac{\dd\mu}{\dd\theta} = 
  v_0 +
  v_1 \mu + 
  v_2 {\mu}^2 . 
\label{eq:Morris_QVF}
\end{equation}

By definition of $\theta \! \left( \mu \right)$ and $\mu \! \left(
\theta \right)$ as inverse functions, it follows that the variance is
also the (geometric and algebraic) inverse of the curvature of the
relative entropy.  We differentiate the second line in
Eq.~(\ref{eq:Lengendre_pair}) twice and substitute
Eq.~(\ref{eq:Morris_QVF}), to produce
\begin{equation}
  \frac{\dd^2 S}{\dd {\mu}^2} = 
  -\frac{\dd\theta}{\dd\mu} = 
  \frac{
    -1
  }{
    v_0 +
    v_1 \mu + 
    v_2 {\mu}^2
  } . 
\label{eq:S_curve_Morris}
\end{equation}

Because we have first and second constants of integration from the
relations following Eq.~(\ref{eq:Lengendre_pair}),
Eq.~(\ref{eq:S_curve_Morris}) has an unambiguous integral.  To assign
meaning to this integral, however, and in the process to expose a
relation between the Morris and Pearson approaches to classification,
we first factor the variance function into an overall normalization
and the roots of the polynomial.  Write 
\begin{equation}
  v_0 +
  v_1 \mu + 
  v_2 {\mu}^2 \equiv
  v_2
  \left( 
    \mu - {\mu}_1
  \right)
  \left( 
    \mu - {\mu}_2
  \right) , 
\label{eq:two_pole_expansion}
\end{equation}
with the solutions
\begin{equation}
  {\mu}_{1,2} = 
  - \frac{v_1}{2 v_2} \mp 
  \sqrt{
    {
      \left(
        \frac{v_1}{2 v_2}
      \right)
    }^2 - 
    \frac{v_0}{v_2}
  } .  
\label{eq:two_pole_solns}
\end{equation}
Then the integral of Eq.~(\ref{eq:S_curve_Morris}) becomes 
\begin{eqnarray}
\lefteqn{
  v_2 S =
} & & 
\nonumber \\
& & 
  \left(
    \frac{
      {\mu}_2 - \mu
    }{
      {\mu}_2 - {\mu}_1
    } 
  \right)
  \log 
  \left(
    \frac{
      {\mu}_2 - \mu
    }{
      {\mu}_2 - {\mu}_0
    } 
  \right) + 
  \left(
    \frac{
      \mu - {\mu}_1 
    }{
      {\mu}_2 - {\mu}_1
    } 
  \right)
  \log 
  \left(
    \frac{
      \mu - {\mu}_1
    }{
      {\mu}_0 - {\mu}_1
    } 
  \right) . 
\nonumber \\
\label{eq:S_gen_soln}
\end{eqnarray}
If we denote by $\varphi \equiv \left( \mu - {\mu}_1 \right) / \left(
{\mu}_2 - {\mu}_1 \right)$, the analytic continuation of a partition
of the unit interval, we may write Eq.~(\ref{eq:S_gen_soln}) as 
\begin{eqnarray}
  v_2 S 
& = &
  \left(
    1 - \varphi 
  \right)
  \log 
  \left(
    \frac{
      1 - \varphi 
    }{
      1 - {\varphi}_0
    }
  \right) + 
  \varphi 
  \log 
  \frac{
    \varphi 
  }{
    {\varphi}_0
  } 
\nonumber \\
& = & 
  D \! 
  \left( 
    \vec{\varphi}
  \right| \! 
  \left|
    {\vec{\varphi}}_0
  \right). 
\label{eq:S_gen_compact}
\end{eqnarray}
In the second line we use $\vec{\varphi}$ to stand for the
``probability distribution'' $\left( \varphi, 1-\varphi \right)$ on
two atoms, and likewise for ${\vec{\varphi}}_0$.  $D \!  \left(
\vec{\varphi} \right| \! \left| {\vec{\varphi}}_0 \right)$ is the
Kullback-Leibler divergence of $\vec{\varphi}$ from the distribution
${\vec{\varphi}}_0$ defined by the equilibrium mean ${\mu}_0$ and the
variance function.  The standard form for the curvature of a
Kullback-Leibler divergence $S$ may be written
\begin{equation}
  v_2 
  {
    \left( {\mu}_2 - {\mu}_1 \right)
  }^2
  \frac{
    \dd^2 S
  }{
    {\dd\mu}^2
  } = 
  \frac{
    1 
  }{
    \varphi 
    \left( 1 - \varphi \right)
  } . 
\label{eq:S_KL_standard_div}
\end{equation}

A slight variation on the formula~(\ref{eq:S_gen_soln}), making use of
forms~(\ref{eq:two_pole_solns}) for the roots, the Legendre transform
relations~(\ref{eq:Lengendre_pair}), and the constants of integration,
reads
\begin{eqnarray}
  2 v_2 \psi \! \left( \theta \right) + 
  v_1 \theta 
& = & 
  \log 
  \left(
    \frac{
      \left(
        {\mu}_2 - \mu
      \right)
      \left(
        \mu - {\mu}_1
      \right)
    }{
      \left(
        {\mu}_2 - {\mu}_0
      \right)
      \left(
        {\mu}_0 - {\mu}_1
      \right)
    }
  \right)
\nonumber \\
& = & 
  \log 
  \left(
    \frac{
      \varphi
      \left( 1 - \varphi \right)
    }{
      {\varphi}_0
      \left( 1 - {\varphi}_0 \right)
    }
  \right) . 
\label{eq:psi_from_variance}
\end{eqnarray}
This integral relation between the cumulant-generating function and
the variance function appears as Eq.~3.7 in
Ref.~\cite{morris82natural}.

\subsubsection{Two fundamental NEF-QVF families, and various limits}

Working in terms of the signs and magnitudes of the coefficients
$v_0$, $v_1$, $v_2$, Morris identifies exactly six inequivalent
natural exponential families with quadratic variance functions.  Three
are continuous (Gaussian, gamma, and hyperbolic-cosecant probability
density functions), and three are discrete (binomial,
negative-binomial, and Poisson probability mass functions), up to
offset and scaling of the natural observation $X$.  We will see here
that, working in terms of the analytic structure of the
entropy~(\ref{eq:S_gen_soln}), and a simple classification of the
roots ${\mu}_{1,2}$, we may identify two main classes, corresponding
to the continuous and discrete distributions, and various limiting
forms of these, which complete Morris's families.

The quantity that distinguishes the continuous from the discrete
NEF-QVF families is the discriminant $d \equiv v_1^2 - 4 v_0 v_2 = 4
v_2^2 {\left( {\mu}_2 - {\mu}_1 \right)}^2$ (which is unchanged by
offset of $X$).  In the case where $d > 0$, the variance
function~(\ref{eq:Morris_QVF}) has two real roots, while if $d < 0$,
it has two complex-conjugate roots.  By choice of offset and scale, we
may obtain Morris's canonical families by making the complex-conjugate
roots purely imaginary when $d < 0$, or by taking one of the two real
roots to lie at the origin if $d > 0$.  

We begin with the imaginary roots, which select the continuous-valued
NEF-QVF distributions.  The canonical form for these is obtained when
$v_1 \equiv 0$, and $v_0, v_2 > 0$.  We may then define 
\begin{equation}
  {\mu}_{1,2} \equiv 
  \mp i \Lambda , 
\label{eq:imag_roots}
\end{equation}
with $\Lambda \equiv \sqrt{v_0 / v_2}$.  

The relative entropy, about a distribution $p_{x \mid 0}$ in the
NEF-QVF family with mean ${\mu}_0$, must have the form
\begin{eqnarray}
  v_2 S 
& = & 
  \frac{1}{2}
  \log 
  \left(
    \frac{
      {\Lambda}^2 + 
      {\mu}^2
    }{
      {\Lambda}^2 + 
      {\mu}_0^2
    }
  \right) + 
  \frac{\mu}{\Lambda}
  \left[
    {\tan}^{-1}
    \left(
      \frac{{\mu}_0}{\Lambda}
    \right) - 
    {\tan}^{-1}
    \left(
      \frac{\mu}{\Lambda}
    \right)
  \right] 
\nonumber \\
& = & 
  \frac{1}{2}
  \log 
  \left(
    \frac{
      {\Lambda}^2 + 
      {\mu}^2
    }{
      {\Lambda}^2 + 
      {\mu}_0^2
    }
  \right) + 
  \frac{\mu}{\Lambda}
  \left[
    {\tan}^{-1}
    \left(
      \frac{\Lambda}{\mu}
    \right) - 
    {\tan}^{-1}
    \left(
      \frac{\Lambda}{{\mu}_0}
    \right) 
  \right] . 
\nonumber \\
\label{eq:S_imag_form}
\end{eqnarray}
The relation of $\theta$ to $\mu$ and ${\mu}_0$ is 
\begin{equation}
  v_2 \theta = 
  \frac{1}{\Lambda}
  \left[
    {\tan}^{-1}
    \left(
      \frac{\mu}{\Lambda}
    \right) - 
    {\tan}^{-1}
    \left(
      \frac{{\mu}_0}{\Lambda}
    \right)  
  \right] . 
\label{eq:theta_imag_form}
\end{equation}
If we choose a background in which ${\mu}_0 = 0$ (by freedom to offset
$X$), it follows that we may write the cumulant-generating function as
\begin{equation}
  v_2 \psi = 
  \frac{1}{2}
  \log 
  \left(
    1 + {\tan}^2
    \left(
      v_2 \Lambda \theta 
    \right)
  \right) . 
\label{eq:psi_imag_form}
\end{equation}

The canonical normalization for this family of distributions is given
by $v_2 = 1$.  One may check directly that they are produced by the
family of hyperbolic-cosecant density functions 
\begin{equation}
  p_{x \mid 0} = 
  \frac{1}{\Lambda} \, 
  \frac{
    1 
  }{
    e^{\pi x / 2 \Lambda} + 
    e^{- \pi x / 2 \Lambda} 
  }
\label{eq:p_0_CSCH}
\end{equation}
(The proof is by contour integral.  Check that
\begin{eqnarray}
  \cos 
  \left( \Lambda \theta \right)
  e^{\psi \left( \theta \right)} 
& = & 
  \frac{1}{\pi}
  \int_0^{\infty}
  \frac{\dd u}{1 + u^2}
  \left(
    {\left( iu \right)}^{\tilde{\theta}} + 
    {\left( -iu \right)}^{\tilde{\theta}}
  \right)  
\nonumber \\
& = & 
  \frac{1}{\pi}
  \int_{-\infty}^{\infty}
  \frac{
    \dd u \, 
    {\left( iu \right)}^{\tilde{\theta}}
  }{
    1 + u^2
  } = 
  1 , 
\label{eq:contour_psi_check}
\end{eqnarray}
with integration variable $u \equiv e^{\pi x / 2 \Lambda}$ and shifted
parameter $\tilde{\theta} \equiv 2 \Lambda \theta / \pi$.  The contour
that avoids branch cuts, in the log-transform to variables $u$, closes in the negative-imaginary half-plane, encircling the pole $u = -i$.)  The distributions at $\Lambda =
1$ are the canonical densities given in
Ref.~\cite{morris82natural}, Eq.~4.2

It is straightforward to check that, as $\Lambda \rightarrow \infty$,
the relative entropy~(\ref{eq:S_imag_form}) reduces to the form 
\begin{equation}
  S \rightarrow 
  - \frac{
    {
      \left( 
        \mu - {\mu}_0
      \right)
    }^2 
  }{
    2v_0
  } , 
\label{eq:S_Gauss_limit}
\end{equation}
for a Gaussian distribution
\begin{equation}
  p_{x \mid 0} = 
  \frac{
    1
  }{
    \sqrt{2 \pi v_0} 
  }
  e^{
    - {
      \left( 
        x - {\mu}_0
      \right)
    }^2 / 
    2 v_0 
  }
\label{eq:p_Gauss_limit}
\end{equation}
with arbitrary mean.  We have used $v_2 {\Lambda}^2 \equiv v_0$ as
$v_2 \rightarrow 0$.  

In the other limit, as $\Lambda \rightarrow 0$, it is convenient to
take $v_2 = 1 / q \equiv 1 / {\mu}_0$, in which
case we recover the relative entropy
\begin{equation}
  S \rightarrow 
  {\mu}_0 - \mu + 
  {\mu}_0
  \log 
  \left( 
    \frac{\mu}{{\mu}_0}
  \right) , 
\label{eq:S_Gamma_limit}
\end{equation}
appropriate to the standard gamma distribution
\begin{equation}
  p_{x \mid 0} = 
  \frac{
    1 
  }{
    \Gamma \left( q \right)
  }
  x^{\left( q-1 \right)}
  e^{-x} . 
\label{eq:p_Gamma_limit}
\end{equation}
Two of the three continuous-valued NEF-QVF families, therefore, are
degenerate limits of the hyperbolic-cosecant distribution, which
represents the generic case.

The discrete-valued families, following when the variance function has
real roots, may be handled in similar fashion.  We choose canonical
forms by offsetting $x$ to set ${\mu}_1 = 0$, and attain this in the
variance function by taking $v_0 \rightarrow 0$.  The canonical scale
for $x$ is then given by taking $v_1 = 1$.

For the discrete distributions, there are two ``interior'' families of
solutions (the binomial and negative binomial), and one limiting
family (the Poisson) that may be reached from either of them.  The
root ${\mu}_2 = - v_1 / v_2$ in all cases.  To obtain the binomial
distribution on $N$ samples with mean ${\mu}_0 = p N$, 
\begin{equation}
  p_{x \mid 0} = 
  \left(
    \begin{array}{c}
      N \\ x 
    \end{array}
  \right)
  p^x 
  {\left( 1-p \right)}^{N-x} , 
\label{eq:p_binom_limit}
\end{equation}
we take ${\mu}_2 = N$, corresponding to $v_2 = -1/N$.  For this
distribution only, the range is finite, $0 \le x \le N$.  The relative
entropy takes the standard form of a Kullback-Leibler divergence
without extending the definition of $\varphi$ by analytic continuation,
\begin{eqnarray}
  S 
& \rightarrow & 
  -N 
  \left\{ 
    \left(
      1 - \frac{\mu}{N}
    \right)
    \log
    \left(
      \frac{
        1 - \mu / N
      }{
        1 - p
      }
    \right) + 
    \frac{\mu}{N}
    \log 
    \left(
      \frac{\mu / N}{p}
    \right)
  \right\} 
\nonumber \\
& = & 
  -N 
  \left\{ 
    \left(
      1 - \frac{\mu}{N}
    \right)
    \log
    \left(
      \frac{
        1 - \mu / N
      }{
        1 - {\mu}_0 / N
      }
    \right) + 
    \frac{\mu}{N}
    \log 
    \left(
      \frac{\mu}{{\mu}_0}
    \right)
  \right\} 
\nonumber \\
& = & 
  - N 
  D \! 
  \left( 
    \vec{\mu} / N
  \right| \! 
  \left|
    \vec{p}
  \right) . 
\label{eq:S_binom_limit}
\end{eqnarray}

The negative binomial distribution is immediately obtained by taking
$N \rightarrow -N$ in the second line of Eq.~(\ref{eq:S_binom_limit})
while holding ${\mu}_0$ fixed.  The corresponding distribution is 
\begin{equation}
  p_{x \mid 0} = 
  \left(
    \begin{array}{c}
      N-1+x \\ x 
    \end{array}
  \right)
  p^x 
  {\left( 1-p \right)}^N , 
\label{eq:p_neg_binom_limit}
\end{equation}
with $p = {\mu}_0 / \left( N + {\mu}_0 \right)$.  This is the other
``interior'' solution, with ${\mu}_2 = -N$ and therefore $v_2 = 1/N$.  

The Poisson distribution is the limit of either of the previous two
forms as $v_2 \rightarrow 0$, so ${\mu}_2 \rightarrow \pm \infty$, at
$p = {\mu}_0$ fixed.  The distribution is 
\begin{equation}
  p_{x \mid 0} = 
  e^{-{\mu}_0}
  \frac{
    {\mu}_0^x 
  }{
    x! 
  } , 
\label{eq:p_Poisson_limit}
\end{equation}
and the entropy becomes 
\begin{equation}
  S \rightarrow 
  \mu - 
  {\mu}_0 - 
  \mu
  \log 
  \left( 
    \frac{\mu}{{\mu}_0}
  \right) , 
\label{eq:S_Poisson_limit}
\end{equation}
For either of the negative binomial or the Poisson, the range of $x$
is unbounded, $x \ge 0$.  

The relative entropy expressions~(\ref{eq:S_Gamma_limit},\ref{eq:S_Poisson_limit}) for the gamma and the Poisson distributions are the same functional form, under exchange of the reference mean ${\mu}_0$ with the distribution mean $\mu$.  Their respective distributions are likewise interchanged under exchange of $x$ with ${\mu}_0$, except that in the gamma case~(\ref{eq:p_Gamma_limit}), a further shift ${\mu}_0 \rightarrow {\mu}_0 - 1$ must be performed as well.  We will return to integer shifts of this form in the next section.

(We note that the association of imaginary roots with
continuous-valued distributions, and of real roots with
discrete-valued distributions, is a defining structural feature of
quantum-mechanical distributions for particles with finite temperature
but continuous time-dependence~\cite{mahan2000many}.  This is one of
many interesting connections to the NEF-QVF families that it will not
be possible to explore in this publication.)

\subsection{Relations to the Pearson system through
  large-deviations formulae}

It is instructive to compare the forms for the entropies of the
distributions in the NEF-QVF families to the logarithms of the
probability densities or mass functions themselves.  By virtue of the
entropy as a large-deviations measure \cite{touchette09the-large}, it and the log-probability will
coincide to leading exponential order for sufficiently sharply peaked
distributions.  

The entropy is defined in the Morris system as a second integral of a
rational function with two poles.  The $\log p_{x \mid 0}$ is defined
in the Pearson system similarly, except that it is a first-integral of
a rational function with two poles~\cite{johnson94continuous}.  The difference between these two
degrees of integration leads to non-coincidence of the two families,
though in many parameter limits they overlap.

We begin by comparing the continuous distributions.  For the Gaussian,
the two functions are identical up to a constant
\begin{eqnarray}
  \log 
  p_{x \mid 0}
& = & 
  - \frac{
    {\left( x - {\mu}_0 \right)}^2
  }{
    2v_0
  } - 
  \frac{1}{2}
  \log 
  \left( 
    2 \pi v_0
  \right)
\nonumber \\
  S 
& = & 
  - \frac{
    {\left( \mu - {\mu}_0 \right)}^2
  }{
    2v_0
  } . 
\label{eq:logp_S_compare_Gauss}
\end{eqnarray}
For the standard gamma with mean ${\mu}_0 = q$,
\begin{eqnarray}
  \log 
  p_{x \mid 0}
& = & 
  q-1 - x + 
  \left( q-1 \right)
  \log 
  \left(
    \frac{x}{q-1}
  \right)
\nonumber \\
  S 
& \approx & 
  q - x + 
  q
  \log 
  \left(
    \frac{x}{q}
  \right) , 
\label{eq:logp_S_compare_Gamma}
\end{eqnarray}
in which the $\approx$ in the second line keeps the first two terms in
Stirling's formula for $\log \Gamma \! \left( q \right)$.  The
functions are identical in form but differ by an offset $q \rightarrow
q-1$.  

The hyperbolic cosecant density shows the least similarity in its
domain of small argument.  However, at small $\Lambda$, where it is
sharply peaked, and at fixed $x$ or $\mu$, the following expansion
becomes informative,
\begin{eqnarray}
  \log 
  p_{x \mid 0}
& = & 
  - \frac{
    \pi \left| x \right|
  }{
    2 \Lambda 
  } - 
  \log \Lambda - 
  \log 
  \left(
    1 + 
    e^{
      - \pi \left| x \right| / \Lambda 
    }
  \right)
\nonumber \\
  S 
& = & 
  - \frac{\mu}{\Lambda}
  {\tan}^{-1}
  \left(
    \frac{\mu}{\Lambda}
  \right) - 
  \log \Lambda + 
  \frac{1}{2}
  \log 
  \left( 
    {\mu}^2 + {\Lambda}^2
  \right) . 
\nonumber \\
\label{eq:logp_S_compare_CSCH}
\end{eqnarray}
For $\mu / \Lambda \gg 1$, ${\tan}^{-1} \left( \mu / \Lambda \right) \rightarrow \mbox{sgn} \!
\left( \mu \right) \pi / 2$, giving the same two
leading terms for $x$ and for $\mu$.  

The discrete distributions behave similarly.  For the binomial,
\begin{eqnarray}
  \log 
  p_{x \mid 0}
& \approx &
  - N  
  D \! 
  \left( 
    \frac{\vec{x}}{N}
  \right| \! 
  \left|
    \frac{{\vec{\mu}}_0}{N}
  \right) 
\nonumber \\
  S 
& = & 
  - N  
  D \! 
  \left( 
    \frac{\vec{\mu}}{N}
  \right| \! 
  \left|
    \frac{{\vec{\mu}}_0}{N}
  \right) , 
\label{eq:logp_S_compare_bin}
\end{eqnarray}
and for the Poisson
\begin{eqnarray}
  \log 
  p_{x \mid 0}
& \approx &
  x - {\mu}_0 - 
  x 
  \log 
  \left(
    \frac{x}{{\mu}_0}
  \right) 
\nonumber \\
  S 
& = & 
  \mu - {\mu}_0 - 
  \mu 
  \log 
  \left(
    \frac{\mu}{{\mu}_0}
  \right) , 
\label{eq:logp_S_compare|Poisson}
\end{eqnarray}
where again $\approx$ stands for the first two terms in Stirling's
formula for factorials.  Within these approximations, the two
functions are identical.  The negative binomial differs by terms at
$\mathcal{O} \! \left( x / N \right)$, but within a similar Stirling
approximation, it may be written 
\begin{eqnarray}
  \log 
  p_{x \mid 0}
& \approx & 
  \left( N + x \right)
  \log 
  \left(
    \frac{
      N + x
    }{
      N + {\mu}_0 
    }
  \right) -
  x
  \log 
  \left(
    \frac{x}{{\mu}_0}
  \right)
\nonumber \\
& & 
  \mbox{} + 
  \left( N + x \right)
  \log 
  \left(
    1 - 
    \frac{1}{N + x}
  \right) - 
  N 
  \log 
  \left(
    1 - 
    \frac{1}{N}
  \right)
\nonumber \\
& & 
  \mbox{} - 
  \log 
  \left(
    1 + 
    \frac{x}{N - 1}
  \right) 
\nonumber \\
& \approx & 
  \left( N + x \right)
  \log 
  \left(
    \frac{
      N + x
    }{
      N + {\mu}_0 
    }
  \right) -
  x
  \log 
  \left(
    \frac{x}{{\mu}_0}
  \right)
\nonumber \\
& & 
  \mbox{} - 
  \mathcal{O} \! 
  \left( 
    \frac{x}{N - 1}
  \right)
\nonumber \\
  S 
& = & 
  \left( N + \mu \right)
  \log 
  \left(
    \frac{
      N + \mu 
    }{
      N + {\mu}_0 
    }
  \right) -
  \mu 
  \log 
  \left(
    \frac{\mu}{{\mu}_0}
  \right) . 
\label{eq:logp_S_compare_neg_bin}
\end{eqnarray}
The leading terms, corresponding to the analytic continuation of the
Kullback-Leibler form, again coincide.  The only differences arise
from shifts $N \rightarrow N-1$ in a subset of terms, similar to the
shift $q \rightarrow q-1$ in Eq.~(\ref{eq:logp_S_compare_Gamma}).  

The equivalence of $\log p_{x \mid 0}$ and $S$ to leading exponential
order is a consequence of the \emph{large-deviations property}~\cite{touchette09the-large} for
these distributions.  The cumulant-generating function is the integral
of the shifted density, 
\begin{equation}
  e^{
    \psi \left( \theta \right)
  } = 
  \int \dd x \, 
  p_{x \mid 0}
  e^{\theta x} . 
\label{eq:psi_genform_def}
\end{equation}
The exponential of the entropy cancels the absolute magnitude of the inserted 
weight factor $e^{\theta x}$ near the
maximum of the shifted distribution, because for sharply peaked
distributions the maximum is near $x \approx \mu$, 
\begin{equation}
  e^{
    S \left( \mu \right)
  } = 
  \int \dd x \, 
  p_{x \mid 0}
  e^{
    \theta \left( \mu \right)
    \left( x - \mu \right) 
  } .  
\label{eq:S_genform_def}
\end{equation}
(This property of the entropy is equivalent to that of functions known
as \emph{effective actions}, as developed in
Ref.~\cite{smith10large-deviation}.)  $S \left( \mu \right)$ is therefore
approximately equal to $p_{x \mid 0}$, evaluated at $x \approx \mu$.
Thus, the Morris restriction to quadratic variance functions implies
that $\log p_{x \mid 0}$, at leading order, will equal the analytic
continuation of a function of Kullback-Leibler form.

%==========================================================

%==========================================================
% Back Matter (References and Notes)
%----------------------------------------------------------
% Style and layout of the references
\bibliographystyle{apa}
\bibliography{entropy}

\end{document}